\documentclass{article}

\usepackage{color}
\usepackage{amsmath, amssymb, amsthm}

\newtheorem{theorem}{Theorem}
 
\newtheorem{prop}{Proposition}

\newtheorem*{proof*}{Proof}

\begin{document}

\title{{\bf On the Average Hitting Times \\of the Directed Wheel}
\vspace{10mm}}

\author{Shunya TAMURA \\
Okegawa City Okegawa West Junior High School \\
Kawataya, Okegawa, 363-0027, Japan \\ 
e-mail: shunya.tamura059@gmail.com 
}

\date{\empty }

\maketitle
\begin{abstract}
In this paper, following the paper ``On the average hitting times of the squares of cycles,'' we provide an explicit formula for the average hitting times of a simple random walk on a directed graph with $N$ vertices, where the graph consists of a cycle with a single absorbing vertex at its center, using elementary methods. Also, we show that the expected hitting time can be expressed in terms of the Fibonacci and Lucas numbers in general.
\end{abstract}

\vspace{10mm}

\begin{small}
\end{small}
\begin{small}
\par\noindent
{\bf Abbr. title}: Hitting times of the directed wheel
\par\noindent
\end{small}

\vspace{10mm}

\section{Introduction \label{sec01}}

A random walk is one of the most straightforward stochastic processes and it's a widely studied model. 
In a random walk on a graph $G$, if a particle is located at a vertex $u$ and $v$ is connected to u by $d(u)$ edges, the probability that the particle moves along the edge $uv$ to the neighboring vertex $v$ is $1/d(u)$. 
Furthermore, when $u$ and $v$ are two vertices on the graph $G$, the expected number of steps for a random walk starting from vertex $u$ to reach vertex $v$ for the first time is called the average hitting times from $u$ to $v$. In a random walk on graph $G$, the average hitting times from vertex $u$ to vertex $v$ is denoted as $h(G;u,v)$. 
The hitting time of a random walk has many useful properties. 
For example, the cover time, which is the expected time for a random walk to visit all vertices of a graph, can be bounded above by a function of the largest hitting time between any two vertices and bounded below by a function of the smallest hitting time between any two vertices \cite{Lovasz}. As a result, there has been extensive research on determining the expected hitting times on various graphs \cite{Lovasz, Sharavas, Ald}. However, finding explicit formulas for the average hitting times from one vertex to any other vertex in a general graph is difficult. 
Recently, explicit formulas have been provided for classes of highly symmetric graphs and similar structures \cite{Etal, TY}. 

In this paper, we calculate the average hitting times for a different model using a similar analytical method as studied in \cite{Etal}. Additionally, we calculate the number of spanning trees, which represent the complexity of the target graph. 
Specifically, we consider a directed cycle graph with $N$ vertices and one absorbing vertex at the center. We call this graph wheel and denote this graph by $W_N^D$. In this case, the average hitting times can be expressed as follows.

\begin{theorem}
The exact formula for the average hitting times of a simple random walk on the graph $W_{N}^{D}$ (where $N \geq 3$) is as follows:

\begin{align*}
h\left( W_{N}^{D} ;0,\ell \right)
=\begin{cases}
\displaystyle\frac{3}{F_{N}}( F_{N} -F_{N-2\ell }) & \text{if} \ N=2n+1 \\ \\
\displaystyle\frac{3}{L_{N}}( L_{N} -L_{N-2\ell }) & \text{if} \ N=2n
\end{cases}
\end{align*}
where $F_{i}$ denotes the $i$-th Fibonacci number and $L_{i}$ denotes the $i$-th Lucas number. 
\end{theorem}

In this formula, for graphs with an odd number of vertices, the expression is in terms of Fibonacci numbers, while for graphs with an even number of vertices, it is in terms of Lucas numbers. 
Our proof relies on combinatorial arguments without requiring any spectral graph theory. 
Furthermore, the formula not only suggests a relation between the average hitting times of a simple random walk on $W_{N}^{D}$ and Fibonacci or Lucas numbers, but this relation is also evident from the intermediate results obtained during the proof.

The rest of this paper is organized as follows. 
In Section \ref{sec02}, we define the Fibonacci and Lucas numbers. 
We introduce essential relations between Fibonacci and Lucas numbers necessary for obtaining our results. 
In Section \ref{sec03}, we present the model we consider, explain how to formulate the matrix equations, and describe our approach to solving the problem. 
In Section \ref{sec04}, we provide proof of the theorem we have obtained, considering cases based on whether the number of vertices in the graph is odd or even. 
In Section \ref{sec05}, we use the Matrix-Tree Theorem, proven by Kirchhoff \cite{Kir}, to calculate the number of spanning trees that represent the complexity of the graph. 
Specifically, we determine the number of inward and outward-spanning trees rooted at a particular vertex in the directed graph $W_N^D$.
Section \ref{sec06} summarizes our results. 

\section{ Fibonacci and Lucas numbers \label{sec02}}

In this section, we first introduce the definitions of the Fibonacci and Lucas sequences. Additionally, we introduce several well-known formulas related to them. The Fibonacci sequence is defined by the recurrence relation $F_{0} =0$, $F_{1} =1$, and $F_{i+2} =F_{i+1} +F_{i}$, where $i=0, 1, 2, \ldots$. The Lucas sequence is similarly defined by the recurrence relation $L_{0} =2$, $L_{1} =1$, and $L_{i+2} =L_{i+1} +L_{i}$, where $i=0, 1, 2, \ldots$. It is well known that Fibonacci numbers and Lucas numbers are closely related. For example, the relationship $L_{i} =F_{i-1} +F_{i+1}$ holds. These related formulas are used throughout this paper. For the proofs of these formulas, please refer to appropriate literature on Fibonacci and Lucas numbers \cite{Thomas, Koshy}.
\begin{align}
F_{2n-2}-3F_{2n}+F_{2n+2}=0 \label{Fib1} \\
L_{2n-2}-3L_{2n}+L_{2n+2}=0  \label{Fib5} \\
F_{n+m}=F_{n}F_{m+1}+F_{m}F_{n-1} \label{Fib2} \\
\displaystyle \sum_{k=1}^{\ell} F_{2k}=F_{2\ell+1}-1\label{Fib3} \\
\displaystyle \sum_{k=1}^{n-\ell} F_{2k-1}=F_{2n-2\ell} \label{Fib4}
\end{align}

\section{Our model \label{sec03}} 

A graph obtained by adding a single absorbing vertex $v_{x}$ to the interior of an $N$-vertex cycle graph is called a wheel graph $W_{N}^{D}$.

The Laplacian matrix $L_{ij}$ of the graph $W_{N}^{D}$ is defined as follows: 
\begin{center}
$L_{ij}=
\begin{cases}
 -1 & \text{if}\ i\ \neq j \ \text{and}\ (v_{i}, v_{j}) \in E \\
 0 & \text{if}\ i \neq j \ \text{and}\ (v_{i}, v_{j}) \notin E \\
\deg_{W_{N}^{D}}( i) & \text{if}\ i = j
\end{cases}$
\end{center}

Next, the row and column corresponding to the absorbing vertex $v_{x}$ are removed from the Laplacian matrix $L$. Furthermore, the first row and first column of the resulting matrix are also removed, yielding a matrix denoted by $L'$. Let $\vec{h}$ be a column vector whose $i$-th element is $h\left(W_{N}^{D}; 0, i\right)$, and let $\vec{1}$ be a column vector of appropriate dimension with all elements equal to 1. In the case of a random walk on $W_{N}^{D}$, if the random walker has not yet reached the absorbing vertex $v_{x}$, they will move to any adjacent vertex with probability 1. Therefore, for $\forall \ \ell \in \mathbb{Z}_{N} \backslash \{x\}$, the following relation holds: \begin{align*}
h\left(W_{N}^{D} ;0,\ell \right)
=\displaystyle\frac{1}{3}\left(1+h\left(W_{N}^{D} ;-1,\ell \right)+1+h\left(W_{N}^{D} ;1,\ell \right) +1+h\left(W_{N}^{D} ;x,\ell \right)\right)
\end{align*}
Here, note that in the wheel graph $W_{N}^{D}$, $h\left( W_{N}^{D} ;x,\ell \right) =0$. Moreover, since the dihedral group $D_{N}$ acts on $W_{N}^{D}$, for $\forall \ k,\ \ell \in \mathbb{Z}_{N}$, the relation $h\left(W_{N}^{D} ;0,\ell \right)=h\left(W_{N}^{D} ;k+\ell ,k\right)=h\left(W_{N}^{D} ;k,k+\ell \right)$ holds. 
Therefore, for $\forall \ \ell \in \mathbb{Z}_{N} \backslash \{x\}$, 
\begin{align*}
-h\left(W_{N}^{D} ;0,\ell -1\right)+3h\left(W_{N}^{D} ;0,\ell \right)-h\left(W_{N}^{D} ;0,\ell +1\right)=3
\end{align*}
is obtained, leading to $L'\vec{h} =3\vec{1}$. By combining this with the fact that for all $\ell \in \mathbb{Z}_{N}$, $h\left(W_{N}^{D} ;0,\ell \right) =h\left( W_{N}^{D} ;0,N-\ell \right)$, the number of variables in the problem can be reduced by half. Let $L'(i, j)$ denote the $(i, j)$-th element of the matrix $L'$. Define $H_{N}$ as an $\displaystyle \lfloor N/2\rfloor \times \lfloor N/2\rfloor $ matrix, where its $(i, j)$-th element is $L'(i, j)$ for $j=N/2\ (N=2n)$. However, the $(i, j)=(n, n-1)$ element is set to $-2$. In all other cases, the element is set to $L'(i, j)+L'(i, N-j)$. Let $\vec{h'}$ be a vector of dimension $\displaystyle \lfloor N/2\rfloor $ corresponding to the elements $h\left( W_{N}^{D} ;0,i\right)$. Then, $H_{N}\vec{h'} =3\vec{1}$ holds. In our problem, it is sufficient to solve this matrix equation.

\section{Proof of Theorem 1 \label{sec04}} 
In this section, we provide proofs for the average hitting times on $W_N^D$ by separating the cases of odd and even numbers of vertices.

\subsection{Odd number of vertices \label{subsec041}}
Coefficient Matrix $H_{N}$ for Odd vertices $W_{N}^{D}$ (When $N=2n+1$): 
The coefficient matrix $H_{N}$ is defined as follows:

\begin{center}
$H_{N}
=\begin{cases}
3 & \text{if}\ ( i=j\ ,\ (i, j) \neq (n, n)) \\
2 & \text{if}\ ( i=j\ ,\ (i, j)=(n, n)) \\
-1 & \text{if}\ ( |i-j|=1) \\
0 & \text{otherwise}
\end{cases}$
\end{center}

We obtain the inverse of the coefficient matrix $H_N$ as follows.
\begin{prop}
The inverse of the coefficient matrix $H_N$ can be expressed as follows: 
\begin{center}
$\displaystyle H_{N}^{-1}
=\frac{1}{F_{N}}
\begin{cases}
F_{N-2i} F_{2j} & \text{if}\ i >j\\
F_{N-2j} F_{2i} & \text{if}\ i\leq j
\end{cases}$
\end{center}
\end{prop}

\noindent
{\bf Proof}. Define the matrix $G_{N}$ as follows:
\begin{align*}
G_{N} =\frac{1}{F_{N}}
\begin{cases}
F_{N-2i} F_{2j} & \text{if}\ i >j\\
F_{N-2j} F_{2i} & \text{if}\ i\leq j 
\end{cases}
\end{align*}

We show that $H_{N} G_{N} =I_{n}$, starting with the diagonal elements being equal to 1.

\begin{itemize}
\item[$\bullet$]
For $(i, j)=(1, 1)=1$:
\begin{align*}
\displaystyle\frac{1}{F_{N}}(3F_{N-2}F_{2}-F_{N-4}F_{2})
&=\frac{1}{F_{N}}(3F_{N-2}-F_{N-4})
=\frac{1}{F_{N}} F_{N}
=1
\end{align*}

\item[$\bullet$]
For $i=j=\ell $ with $(i, j) \notin (1, 1),\ (n, n)$:
\begin{align*}
 & \frac{1}{F_{N}}(-F_{N-2\ell}F_{2\ell-2}+3F_{N-2\ell } F_{2\ell}-F_{N-2-2\ell}F_{2\ell }) \\
=& \frac{1}{F_{N}}(-F_{N-2\ell}F_{2\ell-2}+2F_{N-2\ell}F_{2\ell }+F_{N-1-2\ell}F_{2\ell }) \\
=& \frac{1}{F_{N}}(-F_{N-2\ell}F_{2\ell-2}+F_{N-2\ell}F_{2\ell}-F_{N+1-2\ell}F_{2\ell}) \\ 
=& \frac{1}{F_{N}}( -F_{N-2\ell } F_{2\ell}+F_{N-2\ell }F_{2\ell -1} +F_{N-2\ell } F_{2\ell }+F_{N+1-2\ell } F_{2\ell }) \\
=& \frac{1}{F_{N}}(F_{N-2\ell}F_{2\ell-1}+F_{N+1-2\ell}F_{2\ell })\\
=& \frac{1}{F_{N}} F_{N} 
= 1
\end{align*}
The fifth equality can be derived from the relation of Fibonacci numbers (\ref{Fib2}). 

\item[$\bullet$]
For $(i, j)=(n, n)=1$:
\begin{align*}
\displaystyle\frac{1}{F_{N}}(-F_{1} F_{N-3}+2F_{1}F_{N-1})
&=\displaystyle\frac{1}{F_{N}}(-F_{N-3}+2F_{N-1}) \\
&=\displaystyle\frac{1}{F_{N}}(F_{N-1}-F_{N-3}+F_{N-1}) \\
&=\displaystyle\frac{1}{F_{N}}(F_{N-2}+F_{N-1}) \\
&=\displaystyle\frac{1}{F_{N}} F_{N}
 =1
\end{align*}

Thus, the diagonal elements are all shown to be 1. 
Next, we show that the other elements are 0.

\item[$\bullet$]
For $i=1,\ 1< j\leq n$:
\begin{align*}
3F_{N-2\ell}F_{2}-F_{N-2\ell}F_{4}
&=3F_{N-2\ell}-3F_{N-2\ell }
=0
\end{align*}

\item[$\bullet$]
For $2 \leq i \leq n-1,\ 2 \leq j \leq n$:

When $i<j$:
\begin{align*}
F_{N-2j}(-F_{2i-2}+3F_{2i}-F_{2i+2})=0
\end{align*}

When $i>j$:
\begin{align*}
F_{2j}(-F_{N-2-2i}+3F_{N-2i}-F_{N+2-2i})=0
\end{align*}
Then, by using the Fibonacci identity (\ref{Fib1}), which contains the relevant factors, it is shown to be 0. 

\item[$\bullet$]
For $i=n,\ 1\leq j < n$:
\begin{align*}
F_{3} F_{2j}+2F_{1}F_{2j}
&=-2F_{2j} +2F_{2j}
=0
\end{align*}

This proves that all non-diagonal elements are 0. Therefore, we have shown that $G_{N} =H_{N}^{-1}$. \qed

By combining the proposition and the matrix equation, a general formula for the average hitting times can be obtained.


\noindent
{\bf Proof}. 
\begin{align*}
h_{0, \ell}
=&\frac{3}{F_{N}}( F_{N-2\ell } F_{2} +F_{N-2\ell } F_{4} +F_{N-2\ell } F_{6} +\cdots +F_{N-2\ell } F_{2\ell } \\ 
&\qquad \qquad \qquad +F_{N-2-2\ell } F_{2\ell } +F_{N-4-2\ell } F_{2\ell } +\cdots +F_{5} F_{2\ell } +F_{3} F_{2\ell } +F_{1} F_{2\ell }) \\
=& \frac{3}{F_{N}}\left(\sum_{k=1}^{\ell} F_{N-2\ell}F_{2k}+\sum_{k=1}^{n-\ell}F_{N-2\ell -2k} F_{2\ell }\right) \\
=& \frac{3}{F_{N}}\left( F_{N-2\ell}\sum_{k=1}^{\ell} F_{2k} +F_{2\ell}\sum_{k=1}^{n-\ell}F_{N-2\ell -2k}\right)
\end{align*}

The Fibonacci sequence formula is applied to (\ref{Fib3}) and (\ref{Fib4}), and the calculation is continued.
\begin{align*}
h_{0,\ell }
=& \frac{3}{F_{N}}(F_{N-2\ell}(F_{2\ell+1}-1)+F_{2\ell}F_{N-1-2\ell }) \\
=& \frac{3}{F_{N}}(F_{N-2\ell}F_{2\ell+1}-F_{N-2\ell}+F_{2\ell}F_{N-1-2\ell }) \\
=& \frac{3}{F_{N}}(F_{N}-F_{N-2\ell })
\end{align*}
\end{itemize}

Therefore, we have the desired conclusion. 
\qed

This proof methodically demonstrates the inverse matrix and how it leads to a general formula for the average hitting times in a wheel graph or related structure.

\subsection{Even number of vertices \label{subsec042}}
The coefficient matrix $H_{N}$ for an even vertex wheel graph $W_{N}^{D}$, where $N=2n$, is defined as follows:
\begin{align*}
H_{N}=
\begin{cases}
3 & \text{if} \ (i=j) \\
-2 & \text{if} \ ((i, j)=(n, n-1))\\
-1 & \text{if} \ (|i-j|=1,\ (i, j) \neq (n, n-1))\\
0 & \mbox{otherwise}
\end{cases}
\end{align*}

In this case, the inverse of the coefficient matrix $H_{N}$, denoted as $H_{N}$, is given by the following expression: 
\begin{prop}
The inverse of the coefficient matrix $H_N^{-1}$ can be expressed as follows:
\begin{align*}
H_{N}^{-1}
=\displaystyle\frac{1}{L_{N}}
\begin{cases}
L_{N-2i} F_{2j} & \text{if} \ (i >j) \\
L_{N-2j} F_{2i} & \text{if} \ (i \leq j,\ j\neq n) \\
L_{1} F_{2i} & \text{if} \ (i \leq j,\ j=n)
\end{cases}
\end{align*}
\end{prop}

\noindent
{\bf Proof}. We define the matrix $G'_N$ as follows.
\begin{align*}
G'_{N}=
\displaystyle\frac{1}{L_{N}}
\begin{cases}
L_{N-2i} F_{2j} & \text{if} \ (i >j) \\
L_{N-2j} F_{2i} & \text{if} \ (i \leq j,\ j\neq n) \\
L_{1} F_{2i} & \text{if} \ (i \leq j,\ j=n)
\end{cases}
\end{align*}

To prove that $H_{N}G'_{N}=I_{n}$, we first show that the diagonal elements are equal to 1.
\begin{itemize}
\item[$\bullet$]
For $(i, j)=(1, 1)=1$. 
\begin{align*}
\frac{1}{L_{N}}( 3L_{N-2} F_{2} -L_{N-4} F_{2}) & =\frac{1}{L_{N}}( 3L_{N-2} -L_{N-4})\\
 & =\frac{1}{L_{N}}( 3F_{N-3} +3F_{N-1} -F_{N-5} -F_{N-3})\\
 & =\frac{1}{L_{N}}( 3F_{N-3} -F_{N-5} +3F_{N-1} -F_{2n-3})\\
 & =\frac{1}{L_{N}}( F_{N-1} +F_{N+1})
   =1
\end{align*}

\item[$\bullet$]
For $(i, j)=1$. Let $i=j=\ell$. 
\begin{align*}
 & \frac{1}{L_{N}}( -L_{N-2\ell } F_{2\ell -2} +3L_{N-2\ell } F_{2\ell } -L_{N-2-2\ell } F_{2\ell })\\
= & \frac{1}{L_{N}}( L_{N-2\ell }( F_{2\ell } -F_{2\ell -2}) +L_{N-2\ell } F_{2\ell } +F_{2\ell }( L_{N-2\ell } -L_{N-2-2\ell }))\\
= & \frac{1}{L_{N}}( L_{N-2\ell } F_{2\ell -1} +L_{N-2\ell } F_{2\ell } +F_{2\ell }( F_{N+1-2\ell } +F_{N-1-2\ell } -F_{N-1-2\ell } -F_{N-3-2\ell }))\\
= & \frac{1}{L_{N}}( L_{N-2\ell } F_{2\ell -1} +L_{N-2\ell } F_{2\ell } +F_{2\ell }( F_{N+1-2\ell } -F_{N-3-2\ell }))\\
= & \frac{1}{L_{N}}( L_{N-2\ell } F_{2\ell +1} +F_{2\ell }( F_{N+1-2\ell } -F_{N-3-2\ell }))\\
= & \frac{1}{L_{N}}( F_{2\ell +2} F_{N+1-2\ell } +F_{2\ell +1} F_{N-1-2\ell } -F_{2\ell } F_{N-3-2\ell })\\
= & \frac{1}{L_{N}}( F_{2\ell +2} F_{N+1-2\ell } +F_{2\ell +1} F_{N-1-2\ell } -F_{2\ell } F_{N-1-2\ell } +F_{2\ell } F_{N-2-2\ell })\\
= & \frac{1}{L_{N}}( F_{2\ell +2} F_{N+1-2\ell } +F_{2\ell } F_{N-2-2\ell } +F_{2\ell -1} F_{N-1-2\ell })\\
= & \frac{1}{L_{N}}( F_{2\ell +2} F_{N+1-2\ell } +F_{2\ell } F_{N-2\ell } -F_{2\ell } F_{N-1-2\ell } +F_{2\ell -1} F_{N-1-2\ell })\\
= & \frac{1}{L_{N}}( F_{N+1} +F_{2\ell } F_{N+1-2\ell } -F_{2\ell } F_{N-1-2\ell } +F_{2\ell -1} F_{N-1-2\ell })\\
= & \frac{1}{L_{N}}( F_{N+1} +F_{2\ell } F_{N-2\ell } +F_{2\ell -1} F_{N-1-2\ell })\\
= & \frac{1}{L_{N}}( F_{N+1} +F_{N-1})
= 1
\end{align*}

\item[$\bullet$]
For $(i, j)=(n, n)=1$. 
\begin{align*}
\frac{1}{L_{N}}( -2L_{1} F_{N-2} +3L_{1} F_{N}) & =\frac{1}{L_{N}}( 3F_{N} -F_{N-2} -F_{N-2})\\
 & =\frac{1}{L_{N}}( F_{N+2} -F_{N-2})\\
 & =\frac{1}{L_{N}}( F_{N+1} +F_{N} -F_{N} +F_{N-1})\\
 & =\frac{1}{L_{N}}( F_{N+1} +F_{N-1})
 =1
\end{align*}

Thus, it has been shown that the diagonal elements are equal to 1.
Next, as in the case of odd vertices, we will show that the other elements of the matrix are equal to 0.

\item[$\bullet$]
For $i=1, 1< j \leq n-1$
\begin{align*}
3L_{N-2\ell } F_{2} -L_{N-2\ell } F_{4}
&=3L_{N-2\ell } -3L_{N-2\ell }
=0
\end{align*}
When $i=1,\ j=n$
\begin{align*}
3L_{1} F_{2} -L_{1} F_{4} 
&=3-3
=0
\end{align*}

\item[$\bullet$]
For $2 \leq i \leq n-1, 2 \leq j \leq n-1$

When $i<j$ 
\begin{align*}
L_{N-2j}(-F_{2i-2}+3F_{2i}-F_{2i+2})=0
\end{align*}

When $i >j$
\begin{align*}
F_{2j}( -L_{N-2-2i} +3L_{N-2i} -L_{N+2-2i})=0. 
\end{align*}

Then, by using the Fibonacci identity (\ref{Fib1}) and Lucas identity (\ref{Fib5}), which contains the relevant factors, it is shown to be 0.

\item[$\bullet$]
For $2\leq i\leq n-1, j=n$
\begin{align*}
-L_{1} F_{2\ell -2}+3L_{1}F_{2\ell}-L_{1} F_{2\ell+2}
&=-F_{2\ell -2}+3F_{2\ell}-F_{2\ell+2}
=0
\end{align*}
Therefore, it has been demonstrated that all off-diagonal elements are zero.

\qed 
\end{itemize}

By applying this proposition, a general formula for the average hitting times can be obtained. 

\noindent
{\bf Proof}. 
\begin{align*}
h_{0,\ell } = & \frac{3}{L_{N}}( L_{N-2\ell } F_{2} +L_{N-2\ell } F_{4} +L_{N-2\ell } F_{6} +\cdots +L_{N-2\ell } F_{2\ell }\\
 & \qquad \qquad \qquad +L_{N-2-2\ell } F_{2\ell } +F_{N-4-2\ell } F_{2\ell } +\cdots +L_{4} F_{2\ell } +L_{2} F_{2\ell } +L_{1} F_{2\ell })\\
= & \frac{3}{L_{N}}\left(\sum _{k=1}^{\ell } L_{N-2\ell } F_{2k} +\sum _{k=1}^{n-\ell -1} F_{N-2\ell -2k} F_{2\ell } +F_{2\ell }\right)\\
= & \frac{3}{L_{N}}\left( L_{N-2\ell }\sum _{k=1}^{\ell } F_{2k} +F_{2\ell }\sum _{k=1}^{n-\ell -1} L_{N-2\ell -2k} +F_{2\ell }\right)\\
= & \frac{3}{L_{N}}\left( L_{N-2\ell }\sum _{k=1}^{\ell } F_{2k} +F_{2\ell }\sum _{k=1}^{n-\ell -1} L_{2k} +F_{2\ell }\right)
\end{align*}
The Fibonacci sequence formula is applied to (\ref{Fib3}) and (\ref{Fib4}), and the calculation is continued.
\begin{align*}
h_{0,\ell }
=& \frac{3}{L_{N}}\left(L_{N-2\ell }(F_{2\ell +1} -1)+F_{2\ell }\sum_{k=1}^{n-\ell -1}F_{2k+1}+F_{2\ell }\sum _{k=1}^{n-\ell -1}F_{2k-1} +F_{2\ell }\right)\\
= & \frac{3}{L_{N}}(L_{N-2\ell } F_{2\ell +1} -L_{N-2\ell } +F_{2\ell }( F_{N-2\ell } -1) +F_{2\ell } F_{N-2\ell -2} +F_{2\ell })\\
= & \frac{3}{L_{N}}(F_{N-2\ell +1} F_{2\ell +1}+F_{N-2\ell -1}F_{2\ell+1}-L_{N-2\ell }+F_{2\ell }F_{N-2\ell }+F_{2\ell } F_{N-2\ell -2})\\
= & \frac{3}{L_{N}}((F_{N-2\ell+1}F_{2\ell+1}+F_{2\ell }F_{N-2\ell })+(F_{N-2\ell-1}F_{2\ell+1}+F_{2\ell }F_{N-2\ell -2})-L_{N-2\ell })\\
= & \frac{3}{L_{N}}(F_{N+1}+F_{N-1}-L_{N-2\ell})\\
= & \frac{3}{L_{N}}(L_{N}-L_{N-2\ell })
\end{align*}

Thus, the proof is complete. 
\qed 

By combining Subsection \ref{subsec041} and Subsection \ref{subsec042}, we obtain the following main theorem. 
\begin{align*}
h\left( W_{N}^{D} ;0,\ell \right)
=\begin{cases}
\displaystyle\frac{3}{F_{N}}( F_{N} -F_{N-2\ell }) & \text{if} \ N=2n+1 \\ \\
\displaystyle\frac{3}{L_{N}}( L_{N} -L_{N-2\ell }) & \text{if} \ N=2n
\end{cases}
\end{align*}
where $F_{i}$ denotes the $i$-th Fibonacci number and $L_{i}$ denotes the $i$-th Lucas number.

\section{Spanning tree of $W_N^D$ \label{sec05}} 
In this section, we deal with the spanning trees of $W_{N}^{D}$. 
A spanning tree is a subgraph that includes all the vertices of a given graph and contains no cycles. 
Spanning trees are an important concept in various application fields, such as network design, communication path optimization, and fundamental research in graph theory. 
According to Kirchhoff's Matrix-Tree Theorem, the number of spanning trees in a graph $G$ is given by the determinant of the matrix obtained by removing any one row and column from the Laplacian matrix $L$ \cite{Kir}. 
In the case of the directed wheel $W_{N}^{D}$, the Laplacian matrix is defined as in Section \ref{sec04}. 
However, the weights of the directed edges ($v_{i}$, $v_{j}$) and ($v_{j}$, $v_{i}$) are denoted by $a_{ij}$ and $a_{ji}$, respectively, and these weights are independent. 
If there is no directed edge, the weight is set to 0. 
Furthermore, instead of using the diagonal matrix $D$, we define 
\begin{align*}
\widehat{D} & =(\widehat{D}_{i} \delta _{ij})_{i,j=1}^{N+1} ,\qquad \widehat{D}_{i} =\sum _{k=1}^{N+1} a_{ki}. 
\end{align*}
Then,  another type of Laplacian matrix, $\widehat{L}$, is defined as follows.
\begin{align*}
\widehat{L} &=\widehat{D} -A.\ 
\end{align*}
Here, the matrix $A$ represents the weighted adjacency matrix $A=[ a_{ij}]_{i,j=1}^{N+1}$. 
Additionally, since there are no loops, note that the diagonal elements are zero. 

For example, when the weight of each directed edge is set to 1, the Laplacian matrices $L$ and $\widehat{L}$ of $W_4^D$, which consists of a cycle graph with $4$ vertices and an additional absorbing vertex $v_x$ inside, are as follows.

\begin{align*}
L=
\begin{bmatrix}
3 & -1 & 0 & -1 & -1\\
-1 & 3 & -1 & 0 & -1\\
0 & -1 & 3 & -1 & -1\\
-1 & 0 & -1 & 3 & -1\\
0 & 0 & 0 & 0 & 0
\end{bmatrix} & ,\quad 
\widehat{L}
=\begin{bmatrix}
2 & -1 & 0 & -1 & -1\\
-1 & 2 & -1 & 0 & -1\\
0 & -1 & 2 & -1 & -1\\
-1 & 0 & -1 & 2 & -1\\
0 & 0 & 0 & 0 & 4
\end{bmatrix}
\end{align*}

Note that since the symmetry $a_{ij} \neq a_{ji}$ does not hold, $L \neq \widehat{L}$. 
In the Matrix-Tree Theorem, $L$ and $\widehat{L}$ play different important roles. 

\begin{theorem}[Kirchhoff \cite{Kir}, Tatte \cite{Tatte}]
Let $T_I(v, G)$ denote the set of all inward-directed spanning trees rooted at vertex $v$ in the directed graph $G$, and $T_O(v, G)$ denote the set of all outward-directed spanning trees. 
The weight $w(T)$ of a spanning tree $T$ is defined as $w(T)=\Pi_{e \in T} a(e)$, that is, as the product of the weights $a(e)=a_{ij}$ of each edge $e=(v_i, v_j)$ in $T$. 
In this context, the cofactors $( -1)^{i+j}\mbox{{\bf det}}L^{(i,j)}$ and  $(-1)^{i+j}\mbox{{\bf det}} \widehat{L}^{(i,j)}$ of the Laplacians $L$ and $\widehat{L}$ have the following relations: 
\begin{align*}
T_{I}(v_{i}, G)
=\sum_{T \in T_{I}(v_{i}, G)}w(T), \qquad 
T_{O}(v_{j}, G)
=\sum_{T \in T_{O}(v_{j}, G)}w(T). 
\end{align*}
\end{theorem}

By applying this theorem, we determine the number of directed spanning trees in our model $W_N^D$. 

\begin{prop}
The number of inward-directed spanning trees rooted at vertex $v_x$ in $W_N^D$ can be expressed as follows: 
\begin{align*}
T_{I}(v_{i}, G)
&=\sum_{T \in T_{I}(v_{i}, G)}w(T)
=L_{2N}-2 
\end{align*}
where $L_{i}$ denotes the $i$-th Lucas number. 
\end{prop}

\noindent
{\bf Proof}. 
Remove the row and column corresponding to vertex $v_x$ from the Laplacian matrix $L$. 
The resulting matrix $L'$ corresponds to the Laplacian matrix of the wheel graph $W_N$ in the case of an undirected graph. The following explicit formula for the number of spanning trees $\tau(W_N)$ in $W_N$ (where there are no loops or multiple edges) is provided in \cite{Gabow}. 
\begin{align*}
\tau ( W_{N}) & =L_{2N} -2.\ 
\end{align*}
From this, we can immediately obtain the result of the proposition. 
\qed

\begin{prop}
The number of outward-directed spanning trees rooted at vertex $v_x$ and other vertices in $W_N^D$ can be expressed as follows. 
\begin{align*}
T_{O}(v_{j}, G)
&=\begin{cases}
\displaystyle\sum_{T \in T_{O}(v_{j}, G)}w(T)=0 & if\ (\ v_{j}=v_{x}) \\ \\
\displaystyle\sum_{T \in T_{O}(v_{j}, G)}w(T)=N^{2} & if\ (v_{j} \neq v_{x})
\end{cases} .\ 
\end{align*}

\end{prop}

\noindent
{\bf Proof}. 
Since vertex $v_x$ in $W_N^D$ is an absorbing vertex, its out-degree is 0. 
Therefore, there are no outward-directed spanning trees rooted at vertex $v_x$, so $T_{O}(v_{x}, G)=0$. 

Next, consider the outward-directed spanning trees rooted at vertices other than $v_x$ in $W_N^D$. 
Remove one row and column not corresponding to vertex $v_x$ from the Laplacian matrix $\widehat{L}$. 
Then, perform cofactor expansion on the row and column corresponding to vertex $v_x$ in the resulting matrix $\widehat{L}'$. The matrix obtained after this step corresponds to the Laplacian matrix of the cycle graph $C_N$ in the case of an undirected graph. 

Therefore, by combining these facts, we have $T_{O}(v_{j}, G)=N^2$, where $v_j \neq v_x$. 
\qed

\section{Conclusion \label{sec06}} 
In this paper, we derived a formula for the average hitting times of a simple random walk on a directed wheel graph $W_{N}^{D}$, which is formed by attaching a single absorbing vertex $v_{x}$ to the interior of an $N$-vertex cycle graph, based on the method used in \cite{Etal}. 
We also demonstrated that this formula can be expressed using the general terms of the Fibonacci and Lucas sequences, depending on whether the number of vertices is odd or even. 
The relations with electrical circuits have been clarified in models similar to the one we have dealt with. 
Investigating these relations would be an exciting area of research \cite{Louis}. 
The expected hitting time and cover time for the wheel graph with undirected edges have already been studied \cite{Yujun}. Investigating their relations with Fibonacci numbers and Lucas numbers is also an interesting problem \cite{new}.  

\section*{Acknowledgments}
The authors are grateful to Tsuyoshi Miezaki and Yuuho Tanaka for valuable discussions on this subject. 


\end{document}